\documentclass{amsart}

\usepackage{times}
\usepackage{amssymb}
\usepackage{amsmath}
\usepackage{amsthm}
\usepackage{latexsym}
\usepackage{epsf}
\usepackage{graphicx}
\usepackage{color}
\usepackage{setspace}
\usepackage{xypic}
\usepackage{xy}
\usepackage{textcomp}
\usepackage{tikz}
\usepackage{enumerate}

\newcommand{\qedbox}{\nobreak \ifvmode \relax \else
      \ifdim\lastskip<1.5em \hskip-\lastskip
      \hskip1.5em plus0em minus0.5em \fi \nobreak
      \vrule height0.75em width0.5em depth0.25em\fi}
\def\brc#1{\langle #1 \rangle}
\def\noin{\noindent}

\def\dblbar{|\hspace{-.02in}{|}}
\def\strdiv{|_{_{<}}}

\newenvironment{pf}[1][Proof:]{\begin{trivlist}
\item[\hskip \labelsep {\itshape #1}]}{\end{trivlist}}

\theoremstyle{plain}
\newtheorem{theorem}[equation]{Theorem}
\newtheorem{conjecture}[equation]{Conjecture}
\newtheorem{proposition}[equation]{Proposition}

\newtheorem{corollary}[equation]{Corollary}

\theoremstyle{definition}
\newtheorem{definition}[equation]{Definition}

\newtheorem{example}[equation]{Example}

\newtheorem{fact}[equation]{Fact}

\theoremstyle{remark}
\newtheorem{remark}[equation]{Remark}

\numberwithin{equation}{section}

\begin{document}

\thispagestyle{empty}

\title[Koszul Algebra for Monomial Ideals]{On the Koszul Algebra for Trivariate Monomial Ideals}
\author{Jared L Painter}

\begin{abstract}
  We will describe how we can identify the structure of the Koszul algebra for trivariate monomial ideals from minimal free resolutions.  We use recent work of L. Avramov, where he classifies the behavior of Bass numbers of embedding codepth 3 commutative local rings.  His classification relies on a corresponding classification of their respective Koszul algebras, which is comprised of 5 categories.  Using Avramov's classification of the Koszul algebra, along with their respective Bass series we will learn how to identify the Koszul algebra structure by inspecting the minimal free resolution of the quotient ring.  We give a complete classification of the Koszul algebra for generic monomial ideals and offer several examples.  In addition we will describe a class of ideals with a specific Koszul algebra structure which was previously unknown.
\end{abstract}

\maketitle

\section{Introduction}
\label{intro}
Throughout, we will assume that $S=\Bbbk[x,y,z]$ is a trivariate polynomial ring over a field $\Bbbk$, with homogeneous maximal ideal $\mathfrak{m}=\brc{x,y,z}$, and $I\subseteq \mathfrak{m}^2$ is an $\mathfrak{m}$-primary monomial ideal of $S$.  In this case it is known that the embedding codepth $c=3$.  Our goal will be to describe the homology of the Koszul complex $K=K(x,y,z)$ for $R=S/I$.  It is known that $A=H(K)$ is a graded $\Bbbk$ algebra.  We will use Avramov's classification of the Koszul algebra structure for codepth 3 local rings in \cite{LLA12} to describe the Koszul algebra for $R$.  In \cite{LLA12}, Avramov describes five classes of the Koszul algebra structure for codepth 3 local rings, using their corresponding Bass series.  This is done using invariants of the algebra structure along with closed form expressions for the Bass series.  We will learn how to compute some of these invariants from the minimal free resolution of $R$.  We will also provide a complete classification of the Koszul algebra when $I$ is a generic monomial ideal in Theorem~\ref{GenTorClass} and provide a family of new examples in Theorem~\ref{BClass}.  We will conclude the paper by discussing Conjecture~\ref{CompClass}, which offers a complete description of the Koszul algebra for all Artinian trivariate monomial ideals.

We will primarily explore three of the five classes given by Avramov in \cite{LLA12}.  The classes of interest here are, \textbf{T}, \textbf{B}, and \textbf{H}$(p,q)$.  The rings in the class \textbf{C}(3) are known to be complete intersections, as proved by Assmus in \cite{EA}.  For our rings this is straightforward since $R$ is a complete intersection if and only if $I=\brc{x^a,y^b,z^c}$.  The other class to note is \textbf{G}$(r)$ where $r>1$.  In \cite[3.10]{LLA12} it was conjectured that if $R$ is in \textbf{G}$(r)$ with $r>1$, then $R$ must be Gorenstein.  This was shown to be false in \cite{CV}, and it is still unclear if there are any monomial ideals which fall into this class, though it seems unlikely.

Both Theorem~\ref{pDescription} and Theorem~\ref{rDescription} are our primary tools in computing the algebra structure of $A$.  The results of these theorems rely heavily on when we obtain nonzero entries from $I$ in the matrices of the minimal free resolution of $R$.  In \cite[2.6]{JP} it is determined when columns from $f_2$ have only nonzero entries from $I$, which we will need to prove Theorem~\ref{pDescription}.  For this purpose all minimal free resolutions of $R$ given in this paper will be the \emph{maximal ordered resolution} described in \cite[2.3]{JP}.  

From our assumptions above we know that if $I$ is minimally generated by $n$ monomials the minimal free resolution of $R$ will have the form,
\begin{equation}\label{eq:FreeRes}\mathbb{F}:= \ 0 \longrightarrow S^m \stackrel{f_3}{\longrightarrow} S^{m+n-1} \stackrel{f_2}{\longrightarrow} S^n \stackrel{f_1}{\longrightarrow} S \longrightarrow R \longrightarrow 0.\end{equation}
It is also known that we can describe the first two nonzero Bass numbers, $\mu_R^0$ and $\mu_R^1$ from the minimal free resolution.  Specifically, $\mu_R^0=m$ and $\mu_R^1=m+n-1-\hat{r}$ where $\hat{r}$ is the number of rows in the matrix of $f_3$ that are dependant modulo the ideal $I$.  For those that are interested we will give a more detailed description of the computation for $\mu_R^1$ in Section~\ref{Prelim}.  Using this description for $\mu_R^1$ is key to proving Theorem~\ref{rDescription} and ultimately allows us to prove most of our results in Sections 4 and 5.

One of our more interesting results is the class of new examples found in Theorem~\ref{BClass} for the class \textbf{B}.  The initial examples in this class were found by A. Brown in \cite{AEB} and were all of type 2.  Additional examples are found in \cite{CV} of type 1 and 3.  In Theorem~\ref{BClass} we find a class of examples with type $n-3$, where $n\geq 5$ is the number of generators of $I$.  We also use the results of the proof for ~\ref{BClass} to find a family of non-generic examples in the class \textbf{T} in Theorem~\ref{TnonGenClass}.  This along with several other examples seem to indicate that Conjecture~\ref{CompClass} is true, but more machinery is needed for a complete proof of this.

\section{Preliminaries}
\label{Prelim}

In this section we will provide some background on minimal free resolutions for monomial ideals provided in \cite{JP} and on the Koszul algebra.  In addition we will set some notation for the duration of the paper.  We will also revisit some of the results from \cite{LLA12} and provide the necessary components of these results here for the reader's convenience.  Generally a majority of our notation will follow with that set in \cite{LLA12}, with only minor differences.

We begin by setting some notation and making some blanket assumptions.  We will assume that $I=\brc{m_1,\ldots, m_n}$ is a monomial ideal minimally generated by $n$ monomials.  We denote the least common multiple of monomials $m_1,\ldots,m_r$ by $m_{1\ldots r}$, specifically $m_{ij}=[m_i,m_j].$  Throughout this paper the monomial $m_i$ will be represented by $x^{a_i}y^{b_i}z^{c_i}$.  We say that a monomial $m'$ \emph{strongly divides} a monomial $m$, denoted $m'\dblbar m$, if $m'$ divides $m/x_i$ for all variables $x_i$ dividing $m$.  If a monomial $m'$ strictly divides a monomial $m$ we will write $m'\strdiv m$.

Since many of our proofs require information from \cite{JP} we will state some of these definitions here to preserve notation. One such item is the definition of the second syzygies for $R=S/I$ where $I=\brc{m_1,\ldots, m_n}$ from \cite[2.2]{JP} which are given by,
\begin{equation}\label{eq:SecondSyz}\displaystyle \sigma_{ij}=\frac{m_{ij}}{m_j}e_j-\frac{m_{ij}}{m_i}e_i, \textrm{ for } 1\leq i<j\leq n.\end{equation}
It should be noted that the set $\{\sigma_{ij}\}_{i<j}$ is rarely a minimal generating set for the second syzygies.  We will denote the set of all second syzygies by $Z_2$.  Since we are interested in a minimal generating set we will used the \emph{ordered minimal second syzygies} as they are defined in \cite[2.3]{JP}, which is denoted by $S_2$.  To use this definition an ordering must be defined on the $\sigma_{ij}$'s, to stay consistent we will use the standard dictionary ordering on the indices of each $\sigma_{ij}.$  Namely, $\sigma_{ij}<\sigma_{kl}$ if either $i<k$ or $j<l$ when $i=k$. The following definition is a restatement of \cite[2.3]{JP}.
\pagebreak
\begin{definition}
\label{ResCons}
If $I$ is a monomial ideal, then $\sigma_{ij} \in S_2$ if and only if the following conditions are satisfied,
\begin{enumerate}\setlength{\itemsep}{-0pt}
\item $\sigma_{ij} \in Z_2-\mathfrak{m}Z_2$ and
\item $\sigma_{ij} \not= \sum_{k<l} a_{kl}\sigma_{kl}$, $a_{kl} \in S$, in which $\sigma_{ij} < \sigma_{kl}$ for all $k,l$ such that $a_{kl}$ is a unit.
\end{enumerate}
\end{definition}

We will construct all of our free resolutions in this paper to match this definition.
Since many of our results apply to generic monomial ideals we will define them here.  More about the structure for the minimal free resolution of generic monomial ideals can be found in \cite{JP}.

\begin{definition}
\label{DefGeneric}
A monomial ideal $I = \brc{m_1,\ldots,m_n}$ is \emph{generic} if whenever two distinct minimal generators $m_i$ and $m_j$ have the same positive degree in some variable, there is another minimal generator $m_k$ such that $m_k \dblbar m_{ij}$.
\end{definition}

Much is known about generic monomial ideals in general as they have been studied extensively in \cite{BPS}, \cite{EM}, \cite{MS}, \cite{MSY}, and \cite{JP}.  In \cite[3.10 and 4.6]{JP} a very precise description is given for minimal free resolutions of generic monomial ideals.  Moreover we learn how to determine when we get nonzero entries in the matrix $f_3$ from $I$.  The machinery behind this is that minimal resolutions of generic monomial ideals are regular triangulations, see \cite{BPS}.  This means that each column in $f_3$ will contain exactly three nonzero entries.  It is further shown in \cite[4.1]{JP} that when a column in $f_3$ contains only three nonzero entries, these entries must be pure powers of the variables $x$, $y$, and $z$.  These concepts are key in our classification of the Koszul algebra structure for generic monomial ideals in Theorem~\ref{GenTorClass}.

Many of our calculations requires us to compute minimal generating sets for the graded components of $A=H(K)$.  We denote each of these modules by $A_1$, $A_2$, and $A_3$, where $A_i=H_i(K)$.  Since $K$ is the Koszul complex of $R$ over $\mathfrak{m}$, it is just the deleted resolution of $\mathfrak{m}$ over $R$, where each graded module of the resolution is defined by the exterior algebra.  Thus $K$ has the form,
\begin{equation}\label{eq:Koszul}0\setlength\arraycolsep{0.8mm}\small \xrightarrow{\hspace{.15in}} R \xrightarrow{\varphi_3= \left[ \begin{array}{r} z\\-y\\x\end{array}\right]} R^3 \xrightarrow{\varphi_2= \left[ \begin{array}{rrr} -y&-z&0\\x&0&-z\\0&x&y\end{array}\right]} R^3 \xrightarrow{ \varphi_1=\left[ \begin{array}{rrr}x&y&z\end{array}\right]} R \xrightarrow{\hspace{.15in}} 0\end{equation}

\vspace{2mm}

We will choose the standard basis elements $e_1,e_2,e_3$ of $R^3$ for the graded degree 1 part of $K$, and let $e_1\wedge e_2$, $e_1\wedge e_3$ and $e_2\wedge e_3$ be the ordered generators of the graded degree 2 part. We then have that $e_1\wedge e_2 \wedge e_3$ generates the degree 3 part of $K$.  To simplify this notation we set,\\[-12pt]
$$e_{12}=e_1\wedge e_2, \ \ \ \ e_{13}=e_1\wedge e_3, \ \ \ \ e_{23}=e_2\wedge e_3 \ \ \ \textrm{and} \ \ \ e_{123}=e_1\wedge e_2\wedge e_3.$$

\noin Moreover from the differentials in $K$ we see that,\\[-20pt]

$$\renewcommand{\arraystretch}{1.25}\begin{tabular}{@{\hspace{.25cm}}l@{\hspace{.5cm}}@{\hspace{.25cm}}l@{\hspace{.5cm}}@{\hspace{.25cm}}l@{\hspace{.25cm}}}
$\varphi_1(e_1)=x$ & $\varphi_2(e_{12})=xe_2-ye_1$ & \\
$\varphi_1(e_2)=y$ & $\varphi_2(e_{13})=xe_3-ze_1$ & $\varphi_3(e_{123})=ze_{12}-ye_{13}+xe_{23}.$ \\
$\varphi_1(e_3)=z$ & $\varphi_2(e_{23})=ye_3-ze_2$ & \\
\end{tabular}
$$

\vspace{2mm}

\noin These mappings will be important when we begin finding generators for the graded pieces of $A$.  The following fact is also useful for computing the generators of each $A_i$.

\begin{fact}
\label{RankFact}
If $R=S/I$ has minimal free resolution as given in \eqref{eq:FreeRes}, then 
 $$\textrm{rank}_{\Bbbk}(A_1)=n, \ \ \ \ \ \  \textrm{rank}_{\Bbbk}(A_2)=m+n-1 \ \ \ \ \ \textrm{and} \ \ \ \ \  \textrm{rank}_{\Bbbk}(A_3)=m.$$
\end{fact}

We will ultimately move away from direct computations of generators for $A_1$, $A_2$, and $A_3$ as the paper progresses, as this can be complicated.  We will provide a simple way to compute the generators of $A_1$ from the generators of $I$ in Proposition~\ref{A1Gens}.  Using this particular generating set for $A_1$ will simplify many of our proofs.  We will now define the invariants, as given in \cite[1.1]{LLA12}, which we will use to classify $A$ for the remainder of the paper.

\begin{definition}
\label{DefAInv}
For $R=S/I$ and $A=H(K)$ we define the following for $A$,
\begin{enumerate}

\item $p=\textrm{rank}_{\Bbbk}(A_1^2)$
\item $q=\textrm{rank}_{\Bbbk}(A_1\cdot A_2)$
\item $r=\textrm{rank}_{\Bbbk}(\delta_2)$, where $\delta_2:A_2\longrightarrow \textrm{Hom}_{\Bbbk}(A_1,A_3)$ where $\delta_2(x)(y)=xy$ for $x\in A_2$ and $y\in A_1$
\end{enumerate}
\end{definition}

We will use Theorem~\ref{pDescription} to show that $p$ is precisely the number of distinct columns in $f_2$ with only nonzero entries from $I$.  Also Theorem~\ref{rDescription} shows that $r$ is precisely the number of rows in $f_3$ which are dependent modulo $I$.  Although this might not be that surprising, it is useful to think of $r$ in this manner for monomial ideals.  As shown in \cite{JP}, we know that a row in $f_3$ is dependent module $I$, if that row contains a nonzero pure power entry from $I$, when $I$ is generic.  This is what allows us to provide a complete classification for generic monomial ideals in Theorem~\ref{GenTorClass}.

Theorem~\ref{rDescription} shows us that $r=\hat{r}$ in the formula $\mu_R^1=m+n-1-\hat{r}$. Recall we gave a description for $\hat{r}$ as the number of rows in $f_3$ which are dependant modulo $I$, previously in Section~\ref{intro}.  As noted we will provide background for the computation of $\mu_R^1$ and more notably $\hat{r}$ as we will commonly use this to describe $r$.

Since $R$ is local and Artinian the Betti numbers of $\omega_R = \textrm{Ext}_S^3(R,S)$ (the canonical module of $R$) are equal to the Bass numbers of $R$, see \cite{DJL}.
Using this Fact we are able to obtain information about $\mu_R^1$ from the minimal free resolution of $R$. By applying Hom$(\__,S)$ to the deleted resolution for $R$ in \eqref{eq:FreeRes}. This gives us
\begin{equation}\label{eq:HomFree}\small0 \rightarrow \textrm{Hom}(S,S) \stackrel{f_1^*}{\longrightarrow} \textrm{Hom}(S^n,S) \stackrel{f_2^*}{\longrightarrow} \textrm{Hom}(S^{m+n-1},S) \stackrel{f_3^*}{\longrightarrow} \textrm{Hom}(S^m,S) \rightarrow \omega_R \rightarrow 0,\end{equation}
\normalsize which is a minimal free resolution of $\omega_R$ over $S$, see \cite[3.3.9]{BH}.

Now if we tensor Hom$(\mathbb{F},S)$ with $R$ we will get a free presentation of $\omega_R$ as an $R$-module,
\begin{equation}\label{eq:FreePres}\textrm{Hom}(S^{m+n-1},S)\otimes_S R \xrightarrow{f_3^*\otimes_S R} \textrm{Hom}(S^m,S)\otimes_S R \xrightarrow{\hspace{.05in}} \omega_R \xrightarrow{\hspace{.05in}} 0.\end{equation}

Though this is a free presentation of $\omega_R$ as an $R$-module, it may not be minimal.  We may think of each map $f_i^*$ as the transpose of the matrix for $f_i$, and thus the map $f_3^*\otimes_S R$ is the transpose of the matrix of $f_3$ with entries in $R$.  The map $f_3^*\otimes_S R$ will be minimal if the rows of $f_3$ are algebraically independent mod $I$.  To clarify, we say that the $k^{th}$ row of $f_3$, denoted $r_k$, is algebraically \emph{dependent} mod $I$ if there exists $a_i \in S$ such that $r_k - (a_1r_1+ \cdots + a_{k-1}r_{k-1} + a_{k+1}r_{k+1} + \cdots + a_{m+n-1}r_{m+n-1}) \in IS^m$.  This allows us to find an expression for $\mu_R^1$ as the number of algebraically independent rows of $f_3$ mod $I$, which gives the following.

\begin{remark}
\label{Bass1Rmk}
Let $\hat{r}$ be the number of rows in $f_3$ that are dependent mod $I$, then $\mu_R^1=m+n-1-\hat{r}.$
\end{remark}

One might ask how we determine if a row in $f_3$ is dependant mod $I$.  Generally speaking this can be quite complicated even for monomial ideals.  However, it is simple if $I$ is generic. In \cite[4.7]{JP} it is shown that if $I$ is generic, then a row in $f_3$ is dependant mod $I$ if and only if it contains a nonzero pure power generator of $I$.  There are other instances in this paper when we will still use this concept when $I$ is not generic.  Specifically we will use this in both the proofs for Theorems~\ref{BClass} and ~\ref{TnonGenClass}, but it is also more difficult to prove how many rows in $f_3$ are dependant mod $I$ in both of these theorems.

We will now provide information on four of the five classes of the Koszul algebra for codepth 3 rings given in \cite{LLA12}.  Again we omit \textbf{C}(3) from the list since we already know which monomial ideals fall into this class.  The following Theorem is a collection of information from \cite[2.1 and 3.1]{LLA12}.

\begin{theorem}[Avramov]
\label{LAClass}
  Let $R$ be a local ring of codepth 3, with minimal free resolution as given in \eqref{eq:FreeRes} and let $l=n-1$.  The following table lists the Bass series $I_R^R(t)$ along with values for the previously defined invariants for various classes of the Koszul algebra for $R$.

\begin{table}
[htb!]\centering
\label{TorTab1Ex}\caption{Bass series and invariants for Koszul algebra classes}
\renewcommand{\arraystretch}{1.8}
\small\begin{tabular}{@{\hspace{.25cm}}p{3cm}||@{\hspace{.25cm}}p{6cm}|@{\hspace{.25cm}}c@{\hspace{.25cm}}|@{\hspace{.25cm}}c@{\hspace{.25cm}}|@{\hspace{.25cm}}c@{\hspace{.25cm}}}\hline
\textup{Class} & $I_R^R(t)$ & $p$ & $q$ & $r$\\ \hline\hline
$\textrm{\textbf{\textup{T}}}$ & $\displaystyle \frac{m+lt-2t^2-t^3+t^4}{1-t-lt^2-(m-3)t^3-t^5}$ & $3$ & $0$ & $0$\\
$\textrm{\textbf{\textup{B}}}$ & $\displaystyle \frac{m+(l-2)t-t^2+t^4}{1-t-lt^2-(m-1)t^3+t^4}$ & $1$ & $1$ & $2$\\
$\textrm{\textbf{\textup{G}}}(r)$ & $\displaystyle \frac{m+(l-r)t-(r-1)t^2-t^3+t^4}{1-t-lt^2-nt^3+t^4}$ & $0$ & $1$ & $r$\\
$\textrm{\textbf{\textup{H}}}(0,0)$ & $\displaystyle \frac{m+lt+t^2-t^3}{1-t-lt^2-mt^3}$ & $0$ & $0$ & $0$\\
$\textrm{\textbf{\textup{H}}}(p,q) \ p+q\geq 1$ & $\displaystyle \frac{m+(l-q)t-pt^2-t^3+t^4}{1-t-lt^2-(m-p)t^3+qt^4}$ & $p$ & $q$ & $q$ \\
\end{tabular}
\end{table}
\end{theorem}

We would like to note that $R$ is Gorenstein but not a complete intersection if and only if it is in \textbf{G}$(r)$ with $r=n\geq 5$ and $m=1$.  However there are no known examples of our rings given by monomial ideals in \textbf{G}$(r)$.  Moreover we know even if there were monomial ideals in \textbf{G}$(r)$, none of them would be Gorenstein, since our rings are Gorenstein if and only if they are complete intersections, which implies they are in \textbf{C}(3).  Some examples of non-Gorenstein rings have been found in \textbf{G}$(r)$ in \cite{CV}.  Since it does seem unlikely that any of our rings are in \textbf{G}$(r)$ we will ignore these rings for the majority of the paper.

\begin{remark}
\label{GolodRemark}
Another interesting fact is that $R$ is Golod if and only if it is in \textbf{H}$(0,0)$.  We will provide examples of these in Section~\ref{TorExSec}, and will classify all Golod rings for Artinian trivariate monomial ideals when $I$ is generic.
\end{remark}

\section{Computational Methods}
\label{Computation}
In this section we will describe our primary computational tools we will need to identify which Koszul algebra class $R$ falls into.  Many of our proofs in this section require a direct computational approach, so we will begin this section with an example.  In this example we will compute the generators of $A_1$, $A_2$ and $A_3$ along with the invariants $p$, $q$ and $r$, for a particular ring. We will also look at the minimal free resolution of $R$ in this example to see how this ties into the computation of our invariants for $A$.

\begin{example}
\label{TorAlgEx1}
Let $I = \brc{x^3,x^2y,y^3,z^3,x^2z^2}$ and let $R=S/I$.  We know that we can determine the number of minimal generators we will need for $A_1,A_2,$ and $A_3$ by finding the ranks of the respective free modules in a minimal free resolution of $R$.  Computing the minimal resolution of $R$ we get,

\begin{equation}\label{eq:Ex1Res}\small \setlength\arraycolsep{0.8mm} 0 \xrightarrow{\hspace{.1in}}  S^2 \xrightarrow{\left[ \begin{array}{rr}z^2&0\\-y&0\\0&z^3\\x&-y^2z\\0&x^2\\0&y^3\end{array}\right]} S^6 \xrightarrow{\left[ \begin{array}{rrrrrr}-y&-z^2&0&0&0&0\\x&0&-y^2&-z^2&0&0\\0&0&x^2&0&-z^3&0\\0&0&0&0&y^3&-x^2\\0&x&0&y&0&z\end{array}\right]} S^5 \xrightarrow{\hspace{.15in}} \cdots \end{equation}
\normalsize

We know that $A_1$ will always have the same number of minimal generators as $I$, we also have that $\textrm{rank}_{\Bbbk}(A_2)=6$ and $\textrm{rank}_{\Bbbk}(A_3)=2$.  It is possible to compute the generators for $A_1,A_2,$ and $A_3$ in Macaulay 2 by computing the homology of the Koszul complex over $R$.  Here we will explain how to do this by hand.  As usual we must ensure that the set of generators we compute are independent, but we must also ensure that they in the kernel of the differentials on the Koszul complex $K$.  Also since each $A_i$ represents the homology of $K$ we must ensure that the generators we choose for each $A_i$ are not dependent mod im$(\varphi_{i+1})$ from \ref{eq:Koszul}.

The generators of $A_1$ will be comprised of degree 1 elements in $K$ which have the form $m_je_i$ where $m_j$ is a nonzero monomial in $R$ and $i=1,2,3$.  To ensure that this element is in $\textrm{ker}(\varphi_1)$ we must have that $\varphi_1(m_je_i)= 0$ for each generator we choose.  Since $\varphi_1(m_je_i)=m_j\varphi_1(e_i)=m_j\cdot x_i$ this makes it fairly simple to find a minimal generating set for $A_1$.  Specifically we find that,
$$A_1 = \brc{x^2e_1,xye_1,xze_1,y^2e_2,z^2e_3}.$$

It is easy to see that these are all independent with respect to each other but we must also make sure that they are independent with respect to im$(\varphi_2)$.  The easiest way to do this is to view the generators of $A_1$ along with the generators for im$(\varphi_2)$ as $\Bbbk$-vectors.  This gives,
$$\renewcommand{\arraystretch}{1.25}\begin{tabular}{@{\hspace{.15cm}}c@{\hspace{.15cm}}@{\hspace{.15cm}}c@{\hspace{.15cm}}@{\hspace{.15cm}}c@{\hspace{.15cm}}@{\hspace{.15cm}}c@{\hspace{.15cm}}@{\hspace{.15cm}}c@{\hspace{.35cm}}|@{\hspace{.35cm}}c@{\hspace{.15cm}}@{\hspace{.15cm}}c@{\hspace{.15cm}}@{\hspace{.15cm}}c@{\hspace{.15cm}}}
$A_1$ & & & & & $\textrm{im}(\varphi_2)$ & & \\ \hline
$x^2$ & $xy$ & $xz$ & $0$ & $0$ & $-y$ & $-z$ & $0$\\
$0$ & $0$ & $0$ & $y^2$ & $0$ & $x$ & $0$ & $-z$\\
$0$ & $0$ & $0$ & $0$ & $z^2$ & $0$ & $x$ & $y$\\
\end{tabular}
$$
We can see from this that the chosen minimal generators for $A_1$ do indeed represent a minimal generating set.

Finding a minimal generating set for $A_2$ is generally not as simple.  It requires more effort to verify that an element is in $\textrm{ker}(\varphi_2)$ for these generators.  We do know that each of these generators will be made up of degree 2 elements from $K$ which have the form $m_je_{ij}$ with $1\leq i<j\leq 3$.  We want $m_je_{ij}\in \textrm{ker}(\varphi_2)$, so we need $\varphi_2(m_je_{ij})=m_j\varphi_2(e_{ij})=m_j(-x_ie_i+x_je_j)$ to be zero.  Once we have a set of minimal generators we must also make sure that it is independent mod im$(\varphi_3)$.  We choose the following generators for $A_2$,
$$A_2 = \brc{x^2e_{12},xy^2e_{12},x^2ze_{13},xz^2e_{13},x^2ze_{23},y^2z^2e_{23}}.$$
Again to verify the independence of these generators we will view them as $\Bbbk$-vectors along with the generator for im$(\varphi_3)$.
$$\renewcommand{\arraystretch}{1.25}\begin{tabular}{@{\hspace{.15cm}}c@{\hspace{.15cm}}@{\hspace{.15cm}}c@{\hspace{.15cm}}@{\hspace{.15cm}}c@{\hspace{.15cm}}@{\hspace{.15cm}}c@{\hspace{.15cm}}@{\hspace{.15cm}}c@{\hspace{.15cm}}@{\hspace{.15cm}}c@{\hspace{.25cm}}|@{\hspace{.25cm}}c@{\hspace{.15cm}}}
$A_2$ & & & & & & $\textrm{im}(\varphi_3)$ \\ \hline
$x^2$ & $xy^2$ & $0$ & $0$ & $0$ & $0$ & $z$ \\
$0$ & $0$ & $x^2z$ & $xz^2$ & $0$ & $0$ & $-y$ \\
$0$ & $0$ & $0$ & $0$ & $x^2z$ & $y^2z^2$ & $x$ \\
\end{tabular}
$$

We are less concerned with finding a minimal generating set for $A_3$ at the moment, because it is not needed to compute the numbers $p$ and $q$.  We can now compute these numbers using the multiplication tables for $A_1\cdot A_1$ and $A_1\cdot A_2$.

\begin{table}[htb!]\centering
\label{TorTab1Ex}\caption{Example~\ref{TorAlgEx1} -- $A_1^2$}\renewcommand{\arraystretch}{1.25}\begin{tabular}{|@{\hspace{.15cm}}l@{\hspace{.75cm}}||@{\hspace{.5cm}}c@{\hspace{.5cm}}|@{\hspace{.15cm}}c@{\hspace{.15cm}}|@{\hspace{.15cm}}c@{\hspace{.15cm}}|@{\hspace{.15cm}}c@{\hspace{.15cm}}|@{\hspace{.15cm}}c@{\hspace{.15cm}}|} \hline $A_1^2$&$x^2e_1$&$xye_1$&$xze_1$&$y^2e_2$&$z^2e_3$\\ \hline\hline
$x^2e_1$&0&0&0&0&0\\ \hline $xye_1$&0&0&0&0&$xyz^2e_{13}$\\ \hline $xze_1$&0&0&0&$xy^2ze_{12}$&0\\ \hline $y^2e_2$&0&0&$-xy^2ze_{12}$&0&$y^2z^2e_{23}$\\ \hline $z^2e_3$&0&$-xyz^2e_{13}$&0&$-y^2z^2e_{23}$&0\\ \hline
\end{tabular}
\end{table}

\begin{table}[htb!]\centering
\label{TorTab2Ex}\caption{Example~\ref{TorAlgEx1} -- $A_1\cdot A_2$}\renewcommand{\arraystretch}{1.25}\begin{tabular}{|@{\hspace{.1cm}}l@{\hspace{.5cm}}||@{\hspace{.35cm}}c@{\hspace{.35cm}}|@{\hspace{.35cm}}c@{\hspace{.35cm}}|@{\hspace{.35cm}}c@{\hspace{.35cm}}|@{\hspace{.2cm}}c@{\hspace{.2cm}}|@{\hspace{.2cm}}c@{\hspace{.2cm}}|} \hline $A_1\cdot A_2$&$x^2e_1$&$xye_1$&$xze_1$&$y^2e_2$&$z^2e_3$\\ \hline\hline
$x^2e_{12}$&0&0&0&0&0\\ \hline $xy^2e_{12}$&0&0&0&0&$xy^2z^2e_{123}$\\ \hline $x^2ze_{13}$&0&0&0&0&0\\ \hline $xz^2e_{13}$&0&0&0&$-xy^2z^2e_{123}$&0\\ \hline $x^2ze_{23}$&0&0&0&0&0\\ \hline $y^2z^2e_{23}$&0&0&0&0&0\\ \hline
\end{tabular}
\end{table}

The $\Bbbk$-vector space rank of these will be the number of distinct nonzero independent products that we get up to multiplication by a unit.  From this we may find that we are tempted to say that $p=3$ and $q=1$, but this is not the case.  We have yet to check if the degree 2 elements in $A_1\cdot A_1$ and the degree 3 elements in $A_1\cdot A_2$ are nonzero modulo the incoming maps.  Since the incoming map for $A_3$ is just zero we can conclude that $q=1$.  We must however look at the three nonzero products from $A_1\cdot A_1$ along with im$(\varphi_3)$.  This gives the following,
$$\renewcommand{\arraystretch}{1.25}\begin{tabular}{@{\hspace{.25cm}}c@{\hspace{.25cm}}@{\hspace{.25cm}}c@{\hspace{.25cm}}@{\hspace{.25cm}}c@{\hspace{.5cm}}|@{\hspace{.5cm}}c@{\hspace{.25cm}}}
 & & & $\textrm{im}(\varphi_3)$ \\ \hline
$xy^2z$ & $0$ & $0$ & $z$ \\
$0$ & $xyz^2$ & $0$ & $-y$ \\
$0$ & $0$ & $y^2z^2$ & $x$ \\
\end{tabular}
$$
It is clear that we cannot obtain $y^2z^2e_{23}$ from im$(\varphi_3)$ or any of the other two nonzero elements.  We do find that both $xy^2ze_{12}$ and $xyz^2e_{13}$ are zero mod im$(\varphi_3)$.  To see this, we obtain $xy^2ze_{12}$ from im$(\varphi_3)$ by,
$$xy^2\cdot\left(\begin{array}{c}z\\-y\\x\end{array}\right) = \left(\begin{array}{c}xy^2z\\-xy^3\\x^2y^2\end{array}\right) = \left(\begin{array}{c}xy^2z\\0\\0\end{array}\right)$$
since both $-xy^3$ and $x^2y^2$ are in $I$.  We could preform a similar computation on $xyz^2e_{13}$.  Thus we have that $p=1$.

Looking at the information from Theorem~\ref{LAClass} we now have two possibilities, either $R$ is in \textbf{B} or $R$ is in \textbf{H}$(1,1)$.  Computing $r$ from the definition is not always a simple task, so instead we use what we already know about the lower Bass numbers, $\mu_R^0$ and $\mu_R^1$ along with the Bass series given in Theorem~\ref{LAClass}.

We know that $\mu_R^0=m=2$, and looking at the free resolution for $R$ given in \eqref{eq:Ex1Res} we can see that there are exactly two rows in $f_3$ which are dependant mod $I$.  From Remark~\ref{Bass1Rmk} we have that $\mu_S^1=m+n-1-\hat{r}=6-2=4$.  We will now compare this to the expressions for \textbf{B} and \textbf{H}$(1,1)$ in Theorem~\ref{LAClass}.  For \textbf{H}$(1,1)$,\\[10pt]
$I_R^R(t) \ = \ \displaystyle \sum_{i\geq0}{\mu_R^i t^i} \ =  \ \mu_R^0 + \mu_R^1t + \mu_R^2t^2 + \cdots \ = \ \frac{2+3t-t^2-t^3+t^4}{1-t-4t^2-t^3+t^4}$

\vspace{.15in}

$\Longrightarrow \ (\mu_R^0 + \mu_R^1t + \mu_R^2t^2 + \cdots)(1-t-4t^2-t^3+t^4) = 2+3t-t^2-t^3+t^4$

\vspace{.15in}

$\Longrightarrow \ \mu_R^1 - \mu_R^0 = 3 \ \Longrightarrow \ \mu_R^1 = 5.$\\
Since we know that $\mu_R^1=4$ this formula does not work and we can conclude that $I$ is in \textbf{B}.  However we will verify that the Bass series formula will work for \textbf{B}, which gives,\\[10pt]
$I_R^R(t) = \displaystyle \sum_{i\geq0}{\mu_R^i t^i} = \mu_R^0 + \mu_R^1t + \mu_R^2t^2 + \cdots = \frac{2+2t-t^2+t^4}{1-t-4t^2-t^3+t^4}$

\vspace{.15in}

$\Longrightarrow \ (\mu_R^0 + \mu_R^1t + \mu_R^2t^2 + \cdots)(1-t-4t^2-t^3+t^4) = m+2t-t^2+t^4$

\vspace{.15in}

$\Longrightarrow \ \mu_R^1 - \mu_R^0 = 2 \ \Longrightarrow \ \mu_R^1 = 4,$ \ which is what we wanted.
\end{example}

%%%%%%%%%%%%%%%%%%%%%%%%%%%%%%%%%%%%%%%%%%%%%%%%%%%%%%%%%%%%%%%%%%%%%%%%%%%%%%%%%%%%%%%%%%%%%%%%%%%%%%%%%%%%%%%%%%%%%%%%%%%%%%%%%%%%%%%%%%%%%%%%%%%%%%%%%%%%%%%%%%%%%%%%%%%%%%%%%%%%%%%%%%

With Theorem~\ref{LAClass} in mind there are some cases when $n$ is small in which we can easily determine what the Koszul algebra class is for $I$.  These are outlined in \cite[3.4.2 and 3.4.2]{LLA12}.

\begin{fact}
\label{TorAlgSmalln}
The following hold for any codepth 3 local ring $R$.
\begin{enumerate}[1.]
\item If $n=4$ then $R$ is in one of the following classes:
\begin{enumerate}[(i)]
\item \textbf{H}$(3,2)$ with $m=2$.
\item \textbf{T} with $m\geq 3$.
\item \textbf{H}$(3,0)$ with even $m\geq 4$.
\end{enumerate}
\item If $n\geq 5$, $m=2$, and $p>0$, then $R$ is in one of the following classes:
\begin{enumerate}[(i)]
\item \textbf{B} with odd $n$.
\item \textbf{H}$(1,2)$ with even $n$.
\end{enumerate}
\end{enumerate}
\end{fact}

It is not difficult to see that for an $\mathfrak{m}$-primary monomial ideal minimally generated by 4 monomials we must have that either $m=2$ or $m=3$. It can also be shown that if $m=2$ then the largest number of minimal generators that we may have for $I$ is $n=5$, illustrated in Example~\ref{TorAlgEx1}.  We will see that we can only satisfy 1.(i), 1.(ii), and 2.(i) from Fact~\ref{TorAlgSmalln}.  We will provide examples to all of these in Section~\ref{TorExSec}.  Before we begin computing our examples we will find some more efficient way to compute the invariants $p$ and $r$.

We will now begin to relate the invariants for the Koszul algebra to information that we can obtain from the minimal free resolution of $R$.  This will allow us to give more precise descriptions of the Koszul algebra for monomial ideals.  We begin by showing how we can always find a specific minimal generating set for $A_1$ from the minimal generators of $I$.

\begin{proposition}
\label{A1Gens}
Let $I$ be a monomial ideal with minimal generating set $$\{x^a,y^b,z^c,x^{a_1}y^{b_1}z^{c_1},\ldots,x^{a_{\rho}}y^{b_{\rho}}z^{c_{\rho}}\},$$ where $a_j=0$ for $i\leq j \leq \rho,$ and $a_j>0$ for $j<i.$  Then the set $$\{x^{a-1}(e_1),y^{b-1}(e_2),z^{c-1}(e_3),\{x^{a_j-1}y^{b_j}z^{c_j}(e_1)\}_{j=1}^{i-1},\{y^{b_j-1}z^{c_j}(e_2)\}_{j=i}^{\rho}\}$$
is a minimal generating set for $A_1$.
\end{proposition}
\begin{pf}
We already know that if $I$ is minimally generated by $n$ monomials then $A_1$ is minimally generated by $n$ degree 1 elements from $K$.  The minimal generating set we specified for $I$ is general. We have only specified which minimal mixed generators have only positive degrees on $y$ and $z$, and which ones have positive degrees on $x$ with the condition that $a_j=0 \textrm{ for } i\leq j \leq \rho.$  To show that this is a minimal generating set for $A_1$ we will first show that each of the generators are in the kernel of the differential $\varphi_1$ from $K$, then we will show that this set is independent module im$(\varphi_2)$.

It is clear that all of the chosen generators are in $\textrm{ker}(\varphi_1)$, since for any minimal generator of $I$, $m_j$ with positive degree on $x_i$ we have that $\varphi_1(\frac{m_j}{x_i}e_i)=0$.  We can represent all of the mixed generators as $x^{a_j-1}y^{b_j}z^{c_j}(e_1)$ with $1\leq j < i$ and $y^{b_j-1}z^{c_j}(e_2)$ with $i\leq j\leq \rho$.  To check for independence we will view the generators as $\Bbbk$-vectors,
$$\renewcommand{\arraystretch}{1.25}\begin{tabular}{@{\hspace{.15cm}}c@{\hspace{.15cm}}@{\hspace{.15cm}}c@{\hspace{.15cm}}@{\hspace{.15cm}}c@{\hspace{.15cm}}@{\hspace{.15cm}}c@{\hspace{.15cm}}@{\hspace{.15cm}}c@{\hspace{.25cm}}|@{\hspace{.25cm}}c@{\hspace{.15cm}}@{\hspace{.15cm}}c@{\hspace{.15cm}}@{\hspace{.15cm}}c@{\hspace{.15cm}}}
$A_1$ & & & & & $\textrm{im}(\varphi_2)$ & & \\ \hline
$x^{a-1}$ & $0$ & $0$ & $\{x^{a_j-1}y^{b_j}z^{c_j}\}_{1\leq j<i}$ & $0$ & $-y$ & $-z$ & $0$\\
$0$ & $y^{b-1}$ & $0$ & $0$ & $\{y^{b_j-1}z^{c_j}\}_{i\leq j\leq \rho}$ & $x$ & $0$ & $-z$\\
$0$ & $0$ & $z^{c-1}$ & $0$ & $0$ & $0$ & $x$ & $y$\\
\end{tabular}
$$
It is clear that $x^{a-1}e_1,y^{b-1}(e_2),$ and $z^{c-1}(e_3)$ are independent.  For the other generators it is not possible to get any $x^{a_j-1}y^{b_j}z^{c_j}(e_1)$ from another $x^{a_k-1}y^{b_k}z^{c_k}(e_1)$ with $j\not=k$, because this would imply that $x^{a_j}y^{b_j}z^{c_j}$ was a multiple of $x^{a_k}y^{b_k}z^{c_k}$ which cannot happen because they are both minimal generators.  A similar argument can be made for $y^{b_j-1}z^{c_j}(e_2)$.  We can also see that the only way we could generate $x^{a_j-1}y^{b_j}z^{c_j}(e_1)$ from im$(\varphi_2)$ is with the following multiplications,
$$-x^{a_j-1}y^{b_j-1}z^{c_j}\cdot \setlength\arraycolsep{0mm}\small\left(\begin{array}{c}-y\\x\\0\end{array}\right) = \left(\begin{array}{c}x^{a_j-1}y^{b_j}z^{c_j}\\-x^{a_j}y^{b_j-1}z^{c_j}\\0\end{array}\right)$$ \begin{center}or\end{center} $$-x^{a_j-1}y^{b_j}z^{c_j-1}\cdot \left(\begin{array}{c}-z\\0\\x\end{array}\right) = \left(\begin{array}{c}x^{a_j-1}y^{b_j}z^{c_j}\\0\\-x^{a_j}y^{b_j}z^{c_j-1}\end{array}\right)$$
which could only equal $x^{a_j-1}y^{b_j}z^{c_j}(e_1)$ if both $x^{a_j}y^{b_j-1}z^{c_j}$ and $x^{a_j}y^{b_j}z^{c_j-1}$ were in $I$.  But this contradicts the assumption that $x^{a_j}y^{b_j}z^{c_j}$ is a minimal generator of $I$. Similarly we can only generate $y^{b_j-1}z^{c_j}(e_2)$ from im$(\varphi_2)$ with the following,
$$-y^{b_j-1}z^{c_j-1}\cdot \left(\begin{array}{c}0\\-z\\y\end{array}\right) = \left(\begin{array}{c}0\\-y^{b_j-1}z^{c_j}\\y^{b_j}z^{c_j-1}\end{array}\right)$$ which again cannot be equal to $y^{b_j-1}z^{c_j}e_2$ because $y^{b_j}z^{c_j-1}\notin I$.  Thus we have shown that this is indeed a minimal generating set for $A_1$.\qed
\end{pf}

Using Proposition~\ref{A1Gens} were are able to easily find a minimal generating set for $A_1$ given a minimal generating set for $I$.  While this is nice, it still doesn't tell us anything about our invariants for the Koszul algebra.  We will use this proposition to prove our next result, which will directly relate the computation of $p=\textrm{rank}_{\Bbbk}(A_1^2)$ to the number of distinct $\sigma_{ij}\in S_2$ with only nonzero entries from $I$, as described in \cite[2.6]{JP}.

\begin{theorem}
\label{pDescription}
If $I=\brc{m_1,\ldots, m_n}$ is an $\mathfrak{m}$-primary monomial ideal such that, $$\displaystyle m_k'(e_i) = \frac{m_k}{x_i}(e_i) \textrm{ \ and \ } \displaystyle m_l'(e_j) = \frac{m_l}{x_j}(e_j), 1\leq i<j \leq 3, 1\leq k<l \leq n$$ are minimal generators of $A_1$, then $m_k'\cdot m_l'(e_{ij})\in A_1^2 \subseteq A_2$ is not zero if and only if $\sigma_{kl}\in S_2$ and $\sigma_{kl}$ only has nonzero entries from $I$.
\end{theorem}
\begin{pf}
$(\Longrightarrow)$ \ \ Without loss of generality we may assume that $i=1$ and $j=2$.  Then $m_k'(e_1) = x^{a_k-1}y^{b_k}z^{c_k}(e_1)$, $m_l'(e_2) = x^{a_l}y^{b_l-1}z^{c_l}(e_2)$, and by assumption,
$$m_k'\cdot m_l'(e_{12}) \ = \ x^{a_k+a_l-1}y^{b_k+b_l-1}z^{c_k+c_l}(e_{12}) \not= 0.$$
This means that $m_k'\cdot m_l'\notin I$ and that $m_k'\cdot m_l'(e_{12})$ is not dependent mod im$(\varphi_3)$.  To proceed we will construct the generators $m_k$ and $m_l$ so that the previous statement holds.  To construct $m_k$ and $m_l$ so that $m_k'\cdot m_l'\notin I$ we only need to ensure that there is no minimal generator of $I$ which divides $m_k'\cdot m_l'$.  First we notice that neither $m_k$ or $m_l$ may have positive degrees on all variables.  Suppose that $m_k$ does have positive degrees on all variables, that is $a_k,b_k,c_k>0$.  Then we have that $a_l\leq a_k+a_l-1$, $b_l\leq b_k+b_l-1$, and $c_l\leq c_k+c_l$, which implies that $m_l|(m_k'\cdot m_l')$.  It is also easy to see that we get a similar contradiction when $m_k$ and $m_l$ only have positive degrees on the same two variables.  By assumption $m_k$ must have positive degree on $x$, and $m_l$ must have positive degree on $y$.  From this there are only two options, which we will describe by $(m_k,m_l)=\textrm{gcd}(m_k,m_l)$.  Either $(m_k,m_l)$ has positive degree on only $x$ or only $y$, or $(m_k,m_l)=1$.  We will show that the second option is the only possible option.  Suppose that $(m_k,m_l)=x^{\alpha}$ with $\alpha = \textrm{min}\{a_k,a_l\} >0$.  Then without loss of generality we must have that $m_k=x^{a_k}z^{c_k}$ and $m_l=x^{a_l}y^{b_l}$.  From this we see that $$m_k|m_k'\cdot m_l' = x^{a_k+a_l-1}y^{b_l-1}z^{c_k} \textrm{ since } a_k\leq a_k+a_l-1.$$  This is a contradiction which means that $(m_k,m_l)=1$.  Since there is no minimal generator of $I$ that divides $m_k'\cdot m_l'$.  This implies that $\sigma_{kl}\in S_2$ by \cite[2.4]{JP}, and since $(m_k,m_l)=1$ we have that $\sigma_{kl}$ only has nonzero entries from $I$ by \cite[2.6]{JP}.\\

\noindent $(\Longleftarrow)$ \ \ Suppose $\sigma_{kl}\in S_2$ and $\sigma_{kl}$ only has nonzero entries from $I$.  Then by \cite[2.6]{JP} $(m_k,m_l)=1$ and there is no minimal generator $m_i$ of $I$ such that $m_i\strdiv m_{kl}$, by \cite[2.5]{JP}.  We know from construction that one of these minimal generators is a pure power and the other generator only has positive degrees in the other two variables.  Without loss of generality let $m_k=x^a$ and $m_l=y^{b_l}z^{c_l}$ with either $b_l>0$ or $c_l>0$.  By Proposition~\ref{A1Gens} we have that $m_k'(e_1) = x^{a-1}(e_1)$ and $m_l'(e_2) = y^{b_l-1}z^{c_l}(e_2)$ are both minimal generators of $A_1$.  We only need to show that $m_k'\cdot m_l'(e_{12}) \not= 0$ in $A_1^2$.  Computing this we have,
$$m_k'\cdot m_l'(e_{12}) = x^{a-1}y^{b_l-1}z^{c_l}(e_{12}) = \left(\begin{array}{c}x^{a-1}y^{b_l-1}z^{c_l}\\0\\0\end{array}\right).$$  This element of $A_1^2$ is zero if $x^{a-1}y^{b_l-1}z^{c_l}\in I$ or $m_k'\cdot m_l'(e_{12})$ is dependent modulo im$(\varphi_3)$.  Since there is no minimal generator $m_i$ of $I$ such that $m_i\strdiv m_{kl} = x^ay^{b_l}z^{c_l}$ then $m_k'\cdot m_l'\notin I$ since $m_k'\cdot m_l'\strdiv m_{kl}$.  For sake of contradiction if $m_k'\cdot m_l'(e_{12})$ was dependent mod im$(\varphi_3)$ then the only possible way we could obtain $m_k'\cdot m_l'(e_{12})$ from $\varphi_3$ would be,
$$\left(\begin{array}{c}x^{a-1}y^{b_l-1}z^{c_l}\\0\\0\end{array}\right) \ = \ x^{a-1}y^{b_l-1}z^{c_l-1}\cdot \left(\begin{array}{c}z\\-y\\x\end{array}\right) \ = \ \left(\begin{array}{c}x^{a-1}y^{b_l-1}z^{c_l}\\-x^{a-1}y^{b_l}z^{c_l-1}\\x^ay^{b_l-1}z^{c_l-1}\end{array}\right)$$
$$\hspace{7.9cm} = \ \left(\begin{array}{c}x^{a-1}y^{b_l-1}z^{c_l}\\-x^{a-1}y^{b_l}z^{c_l-1}\\0\end{array}\right).$$
For this to be true we must have that $-x^{a-1}y^{b_l}z^{c_l-1}\in I$, which means there is a minimal generator $m_i$ of $I$ such that $m_i|x^{a-1}y^{b_l}z^{c_l-1}$.  But this would also imply that $m_i|(m_k'\cdot m_l')$ which is a contradiction.   Thus we have shown that $m_k'\cdot m_l'(e_{12})\not=0$.  \qed
\end{pf}
 This Theorem provides some informative corollaries.

\begin{corollary}
\label{NumpCor}
Let $I$ be an $\mathfrak{m}$-primary monomial ideal, then $p=\textrm{\textup{rank}}_{\Bbbk} (A_1^2)$ is precisely then number of distinct $\sigma_{ij}\in S_2$ which have only nonzero entries from $I$.
\end{corollary}

The proof of this follows immediately from Theorem~\ref{pDescription}.  This gives us an easy way to find $p$ from the free resolution of $R$.  We only need to determine how many columns of $f_2$ have only entries from $I$. 

\begin{corollary}
\label{MaxpCor}
Let $I$ be an $\mathfrak{m}$-primary monomial ideal, then $p=\textrm{\textup{rank}}_{\Bbbk} (A_1^2)\leq n-1$.
\end{corollary}

This follows immediately from \cite[2.6]{JP} and Theorem~\ref{pDescription}.

It would be nice if we also had way to describe $q$ from the minimal free resolution of $R$.  In general this description is not as clear.  But we can observe that if $R$ is in either \textbf{H}$(p,q)$ or \textbf{T} we have that $q=r$.  Moreover when $R$ is in \textbf{T} or \textbf{B} we know the values of $p,q,$ and $r$ from Theorem~\ref{LAClass}.  So we provide the following proposition which will follow from Theorem~\ref{LAClass}.

\begin{theorem}
\label{rDescription}
Let $I$ be an $\mathfrak{m}$-primary monomial ideal minimally generated by $n$ monomials, such that $R$ is not Gorenstein. Then $r=\textrm{\textup{rank}}_{\Bbbk}(\delta_2)$ is precisely the number of rows in the matrix of $f_3$ which are dependent \textup{mod} $I$.
\end{theorem}
\begin{pf}
Recall from Remark~\ref{Bass1Rmk} we showed that $\mu_R^1=m+n-1-\hat{r}$ where $\hat{r}$ was the number of rows in $f_3$ which were dependent mod $I$. Since $R$ is not Gorenstein we will ignore the class \textbf{C}(3).  Since we are also ignoring \textbf{G}$(r)$ for this paper we only need to show that this holds for \textbf{T}, \textbf{B}, and \textbf{H}$(p,q)$.  Although this result holds for \textbf{G}$(r)$ as well.  To prove this we only need to show that $r=\hat{r}$ from the expressions for the Bass Series in each respective class.

\noindent Class \textbf{T}: \ From the expression for the Bass series for \textbf{T} we have that,\\[10pt]
$I_R^R(t) \ = \ \displaystyle \sum_{i\geq0}{\mu_R^i t^i} \ =  \ \mu_R^0 + \mu_R^1t + \cdots \ = \ \frac{m+lt-2t^2-t^3+t^4}{1-t-lt^2-(m-3)t^3-t^5}$

\vspace{.15in}

$\Longrightarrow \ (\mu_R^0 + \mu_R^1t + \cdots)(1-t-lt^2-(m-3)t^3-t^5) = m+lt-2t^2-t^3+t^4$

\vspace{.15in}

$\Longrightarrow \ \mu_R^1 - \mu_R^0 = l \ \Longrightarrow \ \mu_R^1 = m+n-1.$\\

\noin So there are now rows in $f_3$ which are dependent mod $I$ and consequently for \textbf{T}, $r=0$.\\

\noindent Class \textbf{B}: \ From the expression for the Bass series for \textbf{B} we have that,\\[10pt]
$I_R^R(t) \ = \ \displaystyle \sum_{i\geq0}{\mu_R^i t^i} \ =  \ \mu_R^0 + \mu_R^1t + \cdots \ = \ \frac{m+(l-2)t-t^2+t^4}{1-t-lt^2-(m-1)t^3+t^4}$

\vspace{.15in}

$\Longrightarrow \ (\mu_R^0 + \mu_R^1t + \cdots)(1-t-lt^2-(m-1)t^3+t^4) = m+(l-2)t-t^2+t^4$

\vspace{.15in}

$\Longrightarrow \ \mu_R^1 - \mu_R^0 = l-2 \ \Longrightarrow \ \mu_R^1 = m+n-1-2.$\\

\noin This says that there are exactly two rows in $f_3$ that are dependent mod $I$ and for \textbf{B}, $r=2$.\\

\noindent Class \textbf{H}$(p,q)$: \ For $p=q=0$ we observe that the terms in the numerator and the denominator of the expression for the Bass series with degrees $\leq 1$ are the same as they were in \textbf{T}.  Thus we can conclude for \textbf{H}$(0,0)$ that there are no rows in $f_3$ that are dependent mod $I$. We also know that $r=0$ in this case. If $p+q\geq 1$ we find that the expression for the Bass series for \textbf{H}$(p,q)$ gives,\\[10pt]
$I_R^R(t) \ = \ \displaystyle \sum_{i\geq0}{\mu_R^i t^i} \ =  \ \mu_R^0 + \mu_R^1t + \cdots \ = \ \frac{m+(l-q)t-pt^2-t^3+t^4}{1-t-lt^2-(m-p)t^3+qt^4}$

\vspace{.15in}

$\Longrightarrow \ (\mu_R^0 + \mu_R^1t + \cdots)(1-t-lt^2-(m-p)t^3+qt^4) = m+(l-q)t-pt^2-t^3+t^4$

\vspace{.15in}

$\Longrightarrow \ \mu_R^1 - \mu_R^0 = l-q \ \Longrightarrow \ \mu_R^1 = m+n-1-q.$\\

\noin Thus for \textbf{H}$(p,q)$ we have that there are $q$ rows in $f_3$ which are dependent mod $I$ and we also know that $r=q$. \qed
\end{pf}

\section{Examples}
\label{TorExSec}

In this section we will list several examples of ideals along with their Koszul algebra classification.  We will rely on the computational methods used from Section~\ref{Computation}.  This will allow us to determine the Koszul algebra classification for various ideals by looking at a minimal free resolution.  We will conclude this section by giving a small class of examples with a very specific structure.

\begin{example}
\label{TorAlgEx2}
For this example we will satisfy the conditions 1.(i),1.(ii), and 2.(i) from Fact~\ref{TorAlgSmalln}.  We have already seen that $I = \brc{x^3,x^2y,y^3,z^3,x^2z^2}$ from Example~\ref{TorAlgEx1} satisfies 2.(i) and is in class \textbf{B}.  So we will consider two general ideals.

\noindent \textbf{1.(i):} \ Let $I=\brc{x^a,y^b,z^c,x^{\alpha}y^{\beta}z^{\gamma}}$ with $a,b,c>0$ and exactly one of the degrees $\alpha,\beta,$ or $\gamma$ is zero.  Since the computation of the free resolution for $R$ here requires that we choose which $\alpha,\beta,$ or $\gamma$ is zero, we will assume that $\alpha=0$.  This gives the following minimal free resolution,
$$\setlength\arraycolsep{.8mm}\scriptsize 0 \xrightarrow{\hspace{.2in}}  S^2 \xrightarrow{\left[ \begin{array}{rr}z^{\gamma}&0\\0&y^{\beta}\\-y^{b-\beta}&-z^{c-\gamma}\\x^a&0\\0&x^a\end{array}\right]} S^5 \xrightarrow{\left[ \begin{array}{rrrrr}-y^b&-z^c&-y^{\beta}z^{\gamma}&0&0\\x^a&0&0&-z^{\gamma}&0\\0&x^a&0&0&-y^{\beta}\\0&0&x^a&y^{b-\beta}&z^{c-\gamma}\end{array}\right]} S^4 \xrightarrow{\hspace{.2in}}\cdots$$
We quickly see that $f_2$ has exactly 3 columns with only nonzero entries from $I$ thus $p=3$ by Corollary~\ref{NumpCor}.  We also se that there are exactly 2 rows in $f_3$ which contain only nonzero entries from $I$ so $r=2$ by Proposition~\ref{rDescription}, since $I$ is generic.  Thus we conclude that $I$ is in \textbf{H}$(3,2)$.

\noindent \textbf{1.(ii):} \ Now let $I=\brc{x^a,y^b,z^c,x^{\alpha}y^{\beta}z^{\gamma}}$ with $a,b,c,\alpha,\beta,\gamma>0$. Then the minimal free resolution for $R$ is,
$$\setlength\arraycolsep{.6mm}\scriptsize 0 \xrightarrow{\hspace{.2in}}  S^3 \xrightarrow{\left[ \begin{array}{rrr}z^{\gamma}&0&0\\0&y^{\beta}&0\\-y^{b-\beta}&-z^{c-\gamma}&0\\0&0&x^{\alpha}\\x^{a-\alpha}&0&-z^{c-\gamma}\\0&x^{a-\alpha}&y^{b-\beta}\end{array}\right]} S^6 \xrightarrow{\left[ \begin{array}{rrrrrr}-y^b&-z^c&-y^{\beta}z^{\gamma}&0&0&0\\x^a&0&0&-z^c&-x^{\alpha}z^{\gamma}&0\\0&x^a&0&y^b&0&-x^{\alpha}y^{\beta}\\0&0&x^{a-\alpha}&0&y^{b-\beta}&z^{c-\gamma}\end{array}\right]} S^4 \xrightarrow{\hspace{.2in}}\cdots$$
Since there are three columns in $f_2$ with only nonzero entries from $I$ we have that $p=3$ by Corollary~\ref{NumpCor}.  This time we find that there are no rows in $f_3$ that are dependent mod $I$ which implies that $r=0$ and in turn implies $q=0$.  Thus $I$ is either in \textbf{T} or \textbf{H}$(3,0)$.  From Fact~\ref{TorAlgSmalln} we know that $I$ must be in \textbf{T}.
\end{example}

\begin{remark}
\label{TRemark}
We may want to note in general how we would determine if $I$ is in \textbf{T} or in \textbf{H}$(3,0)$, when we are faced with this question.  The expressions for the Bass series gives that the zeroth and first Bass numbers yield the same result in both cases.  Since it is difficult to compute the higher Bass numbers in general we need a better method to differentiate the two classes.  The class \textbf{T} is studied in \cite{LLA87}.  We refer to \textbf{T} as the truncated exterior algebra.  In \cite[3.5]{LLA87} we learn that one of the properties of this algebra is that the graded $\Bbbk$-algebra $B$ from \cite[1.3]{LLA12} is generated by three distinct elements in degree 1, and generated by the products of the degree 1 generators in degree 2.  For us this translates into having three distinct minimal generators of $A_1$ so that the products of these degree 1 generators are all nonzero and minimally generate $A_1^2$. We will also have that these products are minimal generators of $A_2$.  In the case of 1.(ii) in Example~\ref{TorAlgEx2} we have that the minimal generators $x^a,y^b$, and $z^c$ of $I$ give us three minimal second syzygies, $\sigma_{12},\sigma_{13}$, and $\sigma_{23}$, which have only nonzero entries from $I$.  Connecting this with the proof of Theorem~\ref{pDescription} we have that the minimal generators $x^{a-1}(e_1),y^{b-1}(e_2),$ and $z^{c-1}(e_3)$ of $A_1$ give us the minimal generators of $A_1^2$, which are $x^{a-1}y^{b-1}(e_{12}),x^{a-1}z^{c-1}(e_{13})$, and $y^{b-1}z^{c-1}(e_{23})$.  This description for \textbf{T} will be important when we give our general classification of \textbf{T} for generic monomial ideals.
\end{remark}

\begin{example}
\label{TorAlgEx3}
Let $I=\brc{x^5,y^5,z^5,y^3z^3,xy^4z^2,xy^2z^4}$ which is generic. Then the maps for $f_2$ and $f_3$ from the minimal free resolution of $R$ is,
$$\setlength\arraycolsep{.28mm}\scriptsize S^6 \xrightarrow{\left[ \begin{array}{rrrrrr}z^2&0&0&0&0&0\\0&y^2&0&0&0&0\\0&0&y&z&0&0\\-y&0&-z&0&0&0\\0&-z&0&-y&0&0\\0&0&0&0&x&0\\-x^4&0&0&0&-z&0\\0&0&0&0&0&x\\0&x^4&0&0&0&-y\\0&0&x^4&0&y&0\\0&0&0&x^4&0&z\end{array}\right]} S^{11} \xrightarrow{\left[ \begin{array}{rrrrrrrrrrr}-y^5&-z^5&-y^3z^3&-y^4z^2&-y^2z^4&0&0&0&0&0&0\\x^5&0&0&0&0&-z^3&-xz^2&0&0&0&0\\0&x^5&0&0&0&0&0&-y^3&-xy^2&0&0\\0&0&x^5&0&0&y^2&0&z^2&0&-xy&-xz\\0&0&0&x^4&0&0&y&0&0&z&0\\0&0&0&0&x^4&0&0&0&z&0&y\end{array}\right]} S^6$$
We notice that there are exactly three columns in $f_2$ with only nonzero entries from $I$, so $p=3$ by Corollary~\ref{NumpCor}.  We also notice that there are no rows in $f_3$ with entries from $I$ so $r=0$ since $I$ is generic, which implies that $q=0$.  Thus we are faced with the same dilemma as in the previous example, either $S$ is in \textbf{T} or in \textbf{H}$(3,0)$.  However the minimal generators of $A_1$ which contribute to $\textrm{rank}_{\Bbbk}(A_1^2)$ are $x^4e_1,y^4e_2,y^2z^2e_2,$ and $z^4e_3$.  This does not agree with the description for \textbf{T} in Remark~\ref{TRemark} so we can conclude that $R$ is in \textbf{H}$(3,0)$.
\end{example}

\begin{example}
\label{TorAlgEx4}
For this example we will list several non-generic ideals along with their respective Koszul algebra structure.  The reader is encouraged to check this using the methods we have introduced in this section.
\begin{table}[!htbp]\centering
\label{Ex4Table}\caption{Example~\ref{TorAlgEx4} -- Some Ideals and their Koszul algebras}\small
\renewcommand{\arraystretch}{1.5}\begin{tabular}{|@{\hspace{.15cm}}r@{\hspace{.15cm}}|@{\hspace{.15cm}}l@{\hspace{.5cm}}|@{\hspace{.25cm}}c@{\hspace{.25cm}}|@{\hspace{.25cm}}c@{\hspace{.25cm}}|@{\hspace{.25cm}}c@{\hspace{.25cm}}|@{\hspace{.25cm}}c@{\hspace{.25cm}}|}\hline & Ideal & $p$ & $q$ & $r$ & Class\\ \hline\hline
1.& $I=\brc{x^a,y^b,z^c,x^3y^2z,x^2y^3z,xyz^3}$ & 3 & 0 & 0 & \textbf{T}\\ \hline
2.& $I=\brc{x^a,y^b,z^c,x^3z,y^3z,xyz^3}$ & 1 & 1 & 2 & \textbf{B}\\ \hline
3.& $I=\brc{x^a,y^b,z^c,x^3y^3,x^3y^2z}$ & 2 & 1 & 1 & \textbf{H}(2,1)\\ \hline
4.& $I=\brc{x^a,y^b,z^c,x^3y^3,x^3z^3,y^3z^3,xyz^4}$ & 0 & 0 & 0 & \textbf{H}(0,0)\\ \hline
5.& $I=\mathfrak{m}^i, \ i\geq 1$ & 0 & 0 & 0 & \textbf{H}(0,0)\\ \hline
\end{tabular}\end{table}

\noindent We will be able to see later from some of these examples how it may be difficult to find general classifications for non-generic monomial ideals, but they do seem to fit a pattern, which we will discuss prior to Conjecture~\ref{CompClass}.
\end{example}

The next Proposition gives us a specific class of examples where we know exactly what the zeroth and first Bass numbers are. These are exactly the generic ideals which give us the maximum number of nonzero entries from $I$ in both $f_2$ and $f_3$ from \cite[2.6 and 3.10]{JP}.

\begin{proposition}
\label{TorExProp}
If $I$ is minimally generated by $n=\rho+3\geq 4$ monomials,\\ $\{x^a,y^b,z^c,y^{b_1}z^{c_1},\ldots,y^{b_{\rho}}z^{c_{\rho}}\}$ with $0<b_i<b_{i+1}<b$, and $0<c_{i+1}<c_i<c$ for $1\leq i < \rho$,, then $R$ is in \textbf{\textup{H}}$(p,q)$ with $p=\mu_R^1=n-1$ and $q=\mu_R^0=n-2$.
\end{proposition}
\begin{pf}
Notice that $I$ is minimally generated by $n\geq 4$ and is generic based upon our assumptions.  The free resolution of $R$ is the same as the resolution described in \cite[3.10]{JP} so that we get exactly $n-2$ pure power entries from $I$ in $f_3$.  Using this along with \cite[2.6]{JP} we can see that we have exactly $n-1$ minimal second syzygies with only nonzero entries from $I$.  Thus buy Theorem~\ref{pDescription} and Proposition~\ref{rDescription} we have that $p=n-1\geq 3$ and $r=n-2\geq 2$.  Thus we have that $R$ must be in \textbf{H}$(p,q)$ with $p=n-1$ and $q=r=n-2$.  It is also easy to see that $\mu_S^1=p$ and $\mu_S^0=q$ in this case.\qed
\end{pf}

\section{Classes for Trivariate Monomial Ideals}
For all of the examples we have seen so far, $p\geq q$.  In general this is not true, see Fact~\ref{TorAlgSmalln}.  Currently there are no known examples where $p<q$ for our monomial ideals, but it is unclear as to how we would prove this in general.  However we can give a positive answer to this for generic monomial ideals.

\begin{proposition}
\label{pbigq}
If $I$ is $\mathfrak{m}$-primary and generic then $p\geq r$, which implies $p\geq q$.
\end{proposition}
\begin{pf}
In \textbf{T}, \textbf{B}, and \textbf{H}$(0,0)$ this is clear, so we only need to show this for \textbf{H}$(p,q)$ with $p+q\geq 1$.  For this case $q=r$ and by Proposition~\ref{rDescription} $r$ is precisely the number of rows in $f_3$ which are dependent mod $I$.  For sake of contradiction, suppose that $r>p$. Then there are $r$ rows in $f_3$ which are dependent mod $I$.  Moreover, we know that since $I$ is generic that there are exactly $r$ entries of the same nonzero pure power generator of $I$ in $r$ rows of $f_3$, by \cite[3.9]{JP}.  Suppose that this entry is $x^a$, then we have $r$ sets of generators of the form $\{x^a,y^{b_1}z^{c_1},y^{b_2}z^{c_2}\}, \ldots ,\{x^a,y^{b_r}z^{c_r},y^{b_{r+1}}z^{c_{r+1}}\}$ which all represent minimal third syzygies for $R$.  However this would imply that we have $r+1$ minimal second syzygies generated by the sets of generators $\{x^a,y^{b_1}z^{c_1}\},\ldots , \{x^a,y^{b_{r+1}}z^{c_{r+1}}\}$.  By \cite[2.6]{JP} each of these second syzygies would have only nonzero entries from $I$ and thus $p\geq  r+1$ by Corollary~\ref{NumpCor}, which is a contradiction. Thus $p\geq r$ and hence $p\geq q$.\qed
\end{pf}

We can now give a complete classification of the Koszul algebra for generic monomial ideals.

\begin{theorem}
\label{GenTorClass}
Let $I$ be an $\mathfrak{m}$-primary generic monomial ideal minimally generated by $n$ monomials, then we have the following Koszul algebra classifications for $R$:
\begin{enumerate}
\item If $I=\brc{x^a,y^b,z^c,x^{a_1}y^{b_1}z^{c_1},\ldots x^{a_{\rho}}y^{b_{\rho}}z^{c_{\rho}}}$ with $a_i,b_i,c_i>0$ for all $1\leq i\leq \rho$ and $\rho>0$, then $R$ is in \textbf{\textup{T}}.
\item If $I=\brc{x^a,y^b,z^c,x^{\alpha_1}y^{\beta_1},x^{\alpha_2}z^{\gamma_1},y^{\beta_2}z^{\gamma_2},x^{a_1}y^{b_1}z^{c_1},\ldots x^{a_{\rho}}y^{b_{\rho}}z^{c_{\rho}}}$, $\alpha_i,\beta_i,\gamma_i>0$, then $R$ is in \textbf{\textup{H}}$(0,0)$ and is Golod if and only if there are no $\sigma_{ij}\in S_2$ with only nonzero entries from $I$.
\item Otherwise $R$ is in \textbf{\textup{H}}$(p,q)$ with $p+q\geq 1$.
\end{enumerate}
\end{theorem}
\begin{pf}
\textbf{Case 1:} Let $m_1=x^a$, $m_2=y^b$, and $m_3=z^c$, then it is clear that $\sigma_{12},\sigma_{13},\sigma_{23}\in S_2$ and have only nonzero entries from $I$.  Applying \cite[2.6]{JP} we have that these are the only $\sigma_{ij}\in S_2$ that satisfy this. Thus $p=3$ by Corollary~\ref{NumpCor}.  Also from our construction of $I$ and the results for generic ideals in \cite{JP} we have that there are no rows in $f_3$ that are dependent mod $I$.  By Proposition~\ref{rDescription} this implies that $r=0$ which implies $q=0$ by Theorem~\ref{LAClass}.  Using Theorem~\ref{LAClass} we have that either $R$ is in \textbf{T} or \textbf{H}$(3,0)$.  However since the products $x^{a-1}y^{b-1}(e_{12})$, $x^{a-1}z^{c-1}(e_{13})$, and $y^{b-1}z^{c-1}(e_{23})$ from $A_1\cdot A_1$ generate $A_1^2$ we must have that $R$ is in \textbf{T} by Remark~\ref{TRemark}.

\noindent \textbf{Case 2:}  We know that $R$ is in \textbf{H}$(0,0)$ if and only if $S$ is Golod from Remark~\ref{GolodRemark}.  Both directions of this proof are immediate consequences of Corollary~\ref{NumpCor} and Proposition~\ref{pbigq}.  We note here that the only reason that we state that $x^{\alpha_1}y^{\beta_1},x^{\alpha_2}z^{\gamma_1},y^{\beta_2}z^{\gamma_2}$, $\alpha_i,\beta_i,\gamma_i>0$ are minimal generators of $I$ is because this is required for $p$ to equal zero but does not necessarily imply $p$ is zero.  The condition that there are no $\sigma_{ij}\in S_2$ with only nonzero entries from $I$ is what implies $p=0$.

\noindent \textbf{Case 3:}  To prove this case we need to show that there are no other constructions for $I$ in \textbf{T}, and then we must show that when $I$ is generic we cannot get anything in \textbf{B}.  We know that we cannot have any generic ideals in \textbf{G}$(r)$ since $p\geq r$ when $I$ is generic by Proposition~\ref{pbigq}. From Example~\ref{TorAlgEx3} we saw that we can have an example in \textbf{H}$(3,0)$ having a minimal second syzygies involving the sets of generators $\{x^5,y^5\}$, $\{x^5,z^5\}$, and $\{x^5,y^3z^3\}$.  We know in general that if we have a minimal second syzygy between generators of the form $\{x^a,y^{\beta}z^{\gamma}\}$, $\beta,\gamma>0$, then we cannot have any minimal second syzygies between minimal generators $\{y^b,x^{\alpha}z^{\gamma'}\}$ or $\{z^c,x^{\alpha'}y^{\beta'}\}$.  To have $p=3$ and $q=0$ we must have two other minimal second syzygies from generators of the form, $\{x^a,y^{\beta_1}z^{\gamma_1}\}$ and $\{x^a,y^{\beta_2}z^{\gamma_2}\}$ with $0\leq \beta_1 <\beta < \beta_2 \leq b$ and $0\leq \gamma_2 < \gamma < \gamma_1 \leq c$.  It is clear that this is not in \textbf{T} from Remark~\ref{TRemark} because this implies that four minimal generators of $A_1$ would give us three distinct nonzero elements from $A_1^2$.  Thus the only constructions we have from \textbf{T} are in Case 1.

If $R$ was in \textbf{B} then we would have that $p=1$ and $r=2$.  But this cannot happen by Proposition~\ref{pbigq} since $p\geq r$ when $I$ is generic.  This implies that $R$ cannot be in \textbf{B} when $I$ is generic and the only possibility is that $R$ is in \textbf{H}$(p,q)$ with $p+q\geq 1$.\qed
\end{pf}

We should note that it is possible to find monomial ideals that are in \textbf{B} as we have already done so in Examples~\ref{TorAlgEx1} and \ref{TorAlgEx4}.  The class \textbf{B} comes from special ideals constructed by A. Brown in \cite{AEB}.  All of the ideals given in \cite{AEB} have $m=2$ and it has previously been unknown if any examples with $m>2$ exist.  The following theorem will give us a class of monomial ideals in \textbf{B}, where $m=n-3\geq 2$.  The construction of this class is just an extension of the ideal from Example~\ref{TorAlgEx1}.

%\begin{remark}
%Recall that in Remark~\ref{GenNumEntries} we conjectured that for a generic monomial ideal generated by $n$ elements that we could not have $n-3$ nonzero entries from $I$ in $f_3$.  If this is true %then it would mean that we will never have $q=n-3$, more specifically we would not be able to find generic monomial ideals in the class \textbf{H}$(p,n-3)$.
%\end{remark}

\begin{theorem}
\label{BClass}
Let $I$ be minimally generated by $n=\rho+3$ monomials,\\ $\{x^a,y^b,z^c,x^{a_1}z^{c'},x^{a_2}y^{b_2}z^{c'}, \ldots,x^{a_{\rho-1}}y^{b_{\rho-1}}z^{c'},y^{b_{\rho}}z^{c'}\}$ with $\rho\geq 2$ and $0<c'<c$.
Then $R$ is in \textbf{\textup{B}} with $m = \rho\geq 2$.
\end{theorem}
\begin{pf}
To prove this we will first construct the second and third syzygies in the free resolution of $R$ so that we may compute $p$ and $r$.  We will then only need to show that $q\not=r$ which will prove our result.

We begin by identifying the minimal second syzygies of $R$.  Notice that $[x^a,x^{a_1}z^{c'}]=x^az^{c'}$ and $[x^{a_1}z^{c'},x^{a_i}y^{b_i}z^{c'}]= x^{a_1}y^{b_i}z^{c'}$ both strictly divide $[x^a,x^{a_i}y^{b_i}z^{c'}]= x^ay^{b_i}z^{c'}$, so there are no minimal second syzygies between the minimal generators $x^a$ and $x^{a_i}y^{b_i}z^{c'}$ for $2\leq i \leq \rho$ by \cite[2.4]{JP}.  Similarly we can see that there are no minimal second syzygies between $y^b$ and $x^{a_i}y^{b_i}z^{c'}$ for $1\leq i \leq \rho-1$.  On the other hand we have that $[z^c,x^{a_i}y^{b_i}z^{c'}]=x^{a_i}y^{b_i}z^c$, and the only minimal generators that divide this are $z^c$ and $x^{a_i}y^{b_i}z^{c'}$, which implies that each second syzygy between $z^c$ and $x^{a_i}y^{b_i}z^{c'}$ is minimal for $1\leq i\leq \rho$ by \cite[2.5]{JP}.  In addition we have that each second syzygy between $x^{a_i}y^{b_i}z^{c'}$ and $x^{a_j}y^{b_j}z^{c'}$ is minimal if and only if $|i-j|=1$.  We also notice that the second syzygy between $\{x^a,y^b\}$ is minimal.  This gives us all of the minimal second syzygies for $R$ which will come from the following sets of minimal generators in order, assuming that $f_1$ is ordered in the manner we wrote the generators for $I$ from above,\\[-2pt]
$$\{x^a,y^b\},\{x^a,x^{a_1}z^{c'}\},\{y^b,y^{b_{\rho}}z^{c'}\},\{z^c,x^{a_i}y^{b_i}z^{c'}\}_{1\leq i\leq \rho},$$\begin{center}and\end{center}$$ \{x^{a_i}y^{b_i}z^{c'},x^{a_{i+1}}y^{b_{i+1}}z^{c'}\}_{1\leq i<\rho}.$$\\[-2pt]
This implies that we have exactly $2\rho+2$ minimal second syzygies.

For the minimal third syzygies we need to find the minimal cycles between the sets of minimal second syzygies.  These cycles come from the following sets of generators of $I$ in the given order,\\[-4pt]
$$\small\{x^a,y^b,\{x^{a_i}y^{b_i}z^{c'}\}_{1\leq i\leq \rho}\}, \{z^c,x^{a_i}y^{b_i}z^{c'},x^{a_{i+1}}y^{b_{i+1}}z^{c'}\}_{1\leq i< \rho}.$$\\[-4pt]
\normalsize Counting these we have exactly $\rho$ minimal third syzygies.

It would be difficult to write down a general form of the free resolution for $R$ here so instead we will analyze what kind of entries we would have from the minimal second and third syzygies we have listed above.  From the list for the second syzygies, $\{x^a,y^b\}$ is the only pair in which $(x^a,y^b)=1$.  This implies we have only one minimal second syzygy in $S_2$ with only nonzero entries from $I$, which implies that $p=1$ by Corollary~\ref{NumpCor}.  For the third syzygies we have that all of the entries from the syzygies given by $\{z^c,x^{a_i}y^{b_i}z^{c'},x^{a_{i+1}}y^{b_{i+1}}z^{c'}\}$ will only have nonzero pure power entries of each of the variables by \cite[4.1]{JP}, but none of these entries will also be in $I$.  The minimal third syzygy given by $\{x^a,y^b,\{x^{a_i}y^{b_i}z^{c'}\}_{1\leq i\leq \rho}\}$ will have precisely the following form up to the sign on the entries,\\[-4pt]
$$\setlength\arraycolsep{.8mm}\left(\begin{array}{ccccccccccc}z^{c'}&y^b&x^a&0& \cdots &0&x^{a-a_1}y^{b-b_2}&x^{a-a_2}y^{b-b_3}& \cdots & x^{a-a_{\rho-2}}y^{b-b_{\rho-1}}& x^{a-a_{\rho-1}}y^{b-b_{\rho}}\end{array}\right)$$\\[-4pt] where $a-a_i<a-a_{i+1}$ and $b-b_{i}>b-b_{i+1}$ for $1\leq i\leq \rho-1.$  We note that this is actually a column in $f_3$ we are writing it as a row vector for clarity.  We can already see that we have the entries $y^b$ and $x^a$ in this third syzygy.  In $f_3$ these will entries will be in rows 2 and 3 from our ordering, which correspond with minimal second syzygies between generators $\{x^a,x^{a_1}z^{c'}\}$ and $\{y^b,y^{b_{\rho}}z^{c'}\}$.  Since these are not involved with any of the other cycles we know that all the other values in rows 2 and 3 of $f_3$ will be zeros.  Also since $a-a_i<a-a_{i+1}$ and $b-b_{i}>b-b_{i+1}$ for $1\leq i\leq \rho-1$ none of the other rows from this third syzygy will be dependent mod $I$.  Thus we can conclude that $r=2$ from Proposition~\ref{rDescription}.

We will now compute minimal generating sets for $A_1$ and $A_2$ and show that $q=\textrm{rank}_{\Bbbk}(A_1\cdot A_2) = 1$.  From Proposition~\ref{A1Gens} we have that,\\[-4pt]
$$A_1=\brc{x^{a-1}e_1,y^{b-1}e_2,z^{c-1}e_3,\{x^{a_i-1}y^{b_i}z^{c'}e_1\}_{i=1}^{\rho-1},y^{b_{\rho}-1}z^{c'}e_2}.$$\\[-4pt]  We will show that,\\[-4pt] $$A_2=\langle x^{a-1}y^{b-1}(e_{12}),\{x^{a_i-1}y^{b_{i+1}-1}z^{c'}(e_{12})\}_{i=1}^{\rho-1},x^{a-1}z^{c'-1}(e_{13}),\hspace{1.5cm}$$$$\hspace{2cm}\{x^{a_i-1}y^{b_i}z^{c-1}(e_{13})\}_{i=1}^{\rho-1},y^{b-1}z^{c'-1}(e_{23}), y^{b_{\rho}-1}z^{c-1}(e_{23})\rangle.$$\\[-4pt]  We notice that we have precisely $2\rho+2$ generators here.  It is relatively simple to verify that each of these generators is in $\textrm{ker}(\varphi_2)$, and that all of these generators are linearly independent with each other.  We need to ensure that each of these generators is independent mod im$(\varphi_3)$.  We first consider the following,\\[-4pt]
$$\renewcommand{\arraystretch}{1.25}\begin{tabular}{@{\hspace{.01cm}}c@{\hspace{.01cm}}@{\hspace{.01cm}}c@{\hspace{.01cm}}@{\hspace{.01cm}}c@{\hspace{.01cm}}@{\hspace{.01cm}}c@{\hspace{.01cm}}@{\hspace{.01cm}}c@{\hspace{.01cm}}@{\hspace{.01cm}}c@{\hspace{.01cm}}|@{\hspace{.01cm}}c@{\hspace{.01cm}}}
 & & & & & & $\textrm{im}(\varphi_3)$ \\ \hline
$x^{a-1}y^{b-1}$ & $0$ & $0$ & $0$ & $x^{a_i-1}y^{b_{i+1}-1}z^{c'}$ & $0$ & $z$ \\
$0$ & $x^{a-1}z^{c'-1}$ & $0$ & $0$ & $0$ & $x^{a_i-1}y^{b_i}z^{c-1}$ & $-y$ \\
$0$ & $0$ & $y^{b-1}z^{c'-1}$ & $y^{b_{\rho}-1}z^{c-1}$ & $0$ & $0$ & $x$ \\
\end{tabular}
$$\\[8pt]
It is clear that the first four columns are independent so we only need to check the last two.  To obtain $x^{a_i-1}y^{b_{i+1}-1}z^{c'}(e_{12})$ from im$(\varphi_3)$ we would have to have\\[2pt]
$$x^{a_i-1}y^{b_{i+1}-1}z^{c'-1}\cdot \left(\begin{array}{c}z\\-y\\x\end{array}\right) \ = \ \left(\begin{array}{c}x^{a_i-1}y^{b_{i+1}-1}z^{c'}\\-x^{a_i-1}y^{b_{i+1}}z^{c'-1}\\x^{a_i}y^{b_{i+1}-1}z^{c'-1}\end{array}\right).$$\\[4pt]  But this cannot equal $x^{a_i-1}y^{b_{i+1}-1}z^{c'}(e_{12})$ because neither $x^{a_i-1}y^{b_{i+1}}z^{c'-1}$ nor\\ $x^{a_i}y^{b_{i+1}-1}z^{c'-1}$ are in $I$.

\noindent To obtain $x^{a_i-1}y^{b_i}z^{c-1}(e_{13})$ from im$(\varphi_3)$ we would have to have,\\[2pt]
$$-x^{a_i-1}y^{b_{i}-1}z^{c-1}\cdot \left(\begin{array}{c}z\\-y\\x\end{array}\right) \ = \ \left(\begin{array}{c}-x^{a_i-1}y^{b_i-1}z^c\\x^{a_i-1}y^{b_{i}}z^{c-1}\\-x^{a_i}y^{b_{i}-1}z^{c-1}\end{array}\right) \ = \ \left(\begin{array}{c}0\\x^{a_i-1}y^{b_{i}}z^{c-1}\\-x^{a_i}y^{b_{i}-1}z^{c-1}\end{array}\right).$$\\[4pt]
But $x^{a_i}y^{b_{i}-1}z^{c-1}\notin I$ so this cannot be equal to $x^{a_i-1}y^{b_i}z^{c-1}(e_{13})$.  Since we have show that each arbitrary element is not dependent mod $I$ we can conclude that this is a minimal generating set for $A_2$.

We will now compute $q$ using the multiplication table for $A_1\cdot A_2$.  In this table we will write only the monomial for the multiplication since every element in $A_1\cdot A_2\subseteq A_3\subseteq R$.
\begin{table}[h]\centering
\label{TorTab2B}\caption{Proof of Theorem~\ref{BClass} -- $A_1\cdot A_2$}\renewcommand{\arraystretch}{1.25}\small\begin{tabular}{|@{\hspace{.05cm}}l@{\hspace{.05cm}}||@{\hspace{.05cm}}c@{\hspace{.05cm}}|@{\hspace{.05cm}}c@{\hspace{.05cm}}|@{\hspace{.05cm}}c@{\hspace{.05cm}}|@{\hspace{.05cm}}c@{\hspace{.05cm}}|@{\hspace{.05cm}}c@{\hspace{.05cm}}|} \hline &$x^{a-1}(e_1)$&$x^{a_i-1}y^{b_i}z^{c'}(e_1)$&$y^{b-1}(e_2)$&$y^{b_{\rho}-1}z^{c'}(e_2)$&$z^{c-1}(e_3)$\\ \hline\hline
$x^{a-1}y^{b-1}(e_{12})$&0&0&0&0&0\\ \hline $x^{a_i-1}y^{b_{i+1}-1}z^{c'}(e_{12})$&0&0&0&0&$0$\\ \hline $x^{a-1}z^{c'-1}(e_{13})$&0&0&$x^{a-1}y^{b-1}z^{c'-1}$&0&0\\ \hline $x^{a_i-1}y^{b_i}z^{c-1}(e_{13})$&0&0&0&$0$&0\\ \hline $y^{b-1}z^{c'-1}(e_{23})$&$x^{a-1}y^{b-1}z^{c'-1}$&0&0&0&0\\ \hline $y^{b_{\rho}-1}z^{c-1}(e_{23})$&0&0&0&0&0\\ \hline
\end{tabular}
\end{table}

\noindent Here we have that $A_1\cdot A_2$ is generated by $x^{a-1}y^{b-1}z^{c'-1}(e_{123})$, which implies that $q=1$.  Thus we have shown that $p=q=1$ and $r=2$, therefore $R$ is in \textbf{B}.\qed
\end{pf}

In Theorem~\ref{GenTorClass} we saw what kind of generic monomial ideals would be in \textbf{T}.  We can also find a class of non-generic monomial ideals in \textbf{T} by removing the two mixed double generators of $I$ that we have in Theorem~\ref{BClass}.  This gives us the following theorem.

\begin{theorem}
\label{TnonGenClass}
Let $I$ be minimally generated by $n=\rho+3\geq 4$ monomials,\\ $\{x^a,y^b,z^c,x^{a_1}y^{b_1}z^{c'}, \ldots, x^{a_{\rho}}y^{b_{\rho}}z^{c'}\}$ with $a_i,b_i,c'>0$ for $1\leq i\leq \rho$.
Then $R$ is in \textbf{\textup{T}} with $m = \rho+2\geq 3$.
\end{theorem}
\begin{pf}
The proof here will be similar to what we did in the previous theorem except we will not have to compute the minimal generating sets for $A_1$ and $A_2$.  Without loss of generality we will assume that $$I=\brc{x^a,y^b,z^c,x^{a_1}y^{b_1}z^{c'},\ldots, x^{a_{\rho}}y^{b_{\rho}}z^{c'}}$$ where $c'>0$, $a_i>a_{i+1}>0$ and $0<b_i<b_{i+1}$ for all $1\leq i < \rho$.  Using similar arguments as in the previous theorem we have that the following pairs of minimal generators represent all of the minimal second syzygies for $R$,
$$\{x^a,y^b\},\{x^a,z^c\},\{x^a,x^{a_1}y^{b_1}z^{c'}\},\{y^b,z^c\},\{y^b,x^{a_{\rho}}y^{b_{\rho}}z^{c'}\},\{z^c,x^{a_i}y^{b_i}z^{c'}\}_{i=1}^{\rho},$$\begin{center}and\end{center}$$\{x^{a_i}y^{b_i}z^{c'},x^{a_{i+1}}y^{b_{i+1}}z^{c'}\}_{i=1}^{\rho-1}.$$
We notice that the pairs $\{x^a,y^b\},\{x^a,z^c\},$ and $\{y^b,z^c\}$ represent the only minimal second syzygies which have only nonzero entries from $I$.  Thus $p=3$ by Corollary~\ref{NumpCor}.  We will also note that we only need to show that $q=0$ because this would imply that $R$ satisfies the structure for \textbf{T} by Remark~\ref{TRemark}.

Using the minimal second syzygies from above we find that the minimal cycles for the minimal third syzygies are given by the sets of generators,
\small$$\{x^a,y^b,\{x^{a_i}y^{b_i}z^{c'}\}_{i=1}^{\rho}\},\{x^a,z^c,x^{a_1}y^{b_1}z^{c'}\},\{y^b,z^c,x^{a_{\rho}}y^{b_{\rho}}z^{c'}\},$$\begin{center}and\end{center}$$\{z^c,x^{a_i}y^{b_i}z^{c'},x^{a_{i+1}}y^{b_{i+1}}z^{c'}\}_{i=1}^{\rho-1}.$$
\normalsize All but the first cycle only involve three minimal generators none of which satisfy the criterion needed to admit a nonzero entry from $I$ in $f_3$.  The cycle given by\\ $\{x^a,y^b,\{x^{a_i}y^{b_i}z^{c'}\}_{i=1}^{\rho}\}$ will be given up to a sign on the nonzero entries by,
$$\setlength\arraycolsep{.8mm}\left(\begin{array}{cccccccccccc}z^{c'}&0&y^{b-b_1}&0&x^{a-a_{\rho}}&0& \cdots &0&x^{a-a_1}y^{b-b_2}&x^{a-a_2}y^{b-b_3}& \cdots & x^{a-a_{\rho-1}}y^{b-b_{\rho}}\end{array}\right)$$ where $a-a_i<a-a_{i+1}$ and $b-b_{i}>b-b_{i+1}$ for $1\leq i\leq \rho-1.$  It is clear that none of these nonzero elements are in $I$, and none of them are dependent mod $I$ after row operations.  Thus we have shown that there are no rows in $f_3$ which are dependent mod $I$ which implies that $q=0$. \qed
\end{pf}

In both of the previous theorems the classes are less than robust because generally we would not have ideals generated like this by randomly picking minimal generators.  We would eventually like to have a complete classification of the Koszul algebra for trivariate monomial ideals.  The issue this that in general we do not know how the permissable row operations will affect $r$ for a general monomial ideal.  It is also unclear as to whether or not we can get other nonzero entries from $I$ in $f_3$ for these ideals in general.  These ideas seem to be the key to showing the Koszul Algebra classification in general.

Another interesting observation we note from our class in \textbf{B} is that in $f_3$ we will have exactly two different nonzero pure power entries in the same column, which are also in $I$.  In some sense this explains the reason that $q\not= r$ in \textbf{B}.  One could go back to Example~\ref{TorAlgEx1} and see that we have both $y^3$ and $z^3$ as entries in the same column of $f_3$, and these are both minimal generators of $I$.  Looking at the multiplication between $A_1\cdot A_2$ we can see that both $y^2e_2\cdot xz^2e_{13}$ and $z^2e_3\cdot xy^2e_{12}$ yield a nonzero value in $A_1\cdot A_2$.  However since we get that both of these multiplications are the same $q=1$.  In any case it seems that this situation where $q=1$ and $r=2$ is unique to the ideals having exactly two different pure power entries from $I$ in the same column of $f_3$.

We can also see from Example~\ref{TorAlgEx4} that we can find instances of non-generic monomial ideals in each class.  The pattern does not seem that far off from our generic classification, the only difference being that we have to take the class \textbf{B} into account.  In Example~\ref{TorAlgEx4} we even provide an example in \textbf{B} that is outside the scope of our classification from Theorem~\ref{BClass}.  The common trait here is that in both cases we have exactly two different nonzero pure power entries from $I$ in the same column in $f_3$.  

Upon inspection of any cycle of $R$, with corresponding generators $\{m_1,\ldots, m_r\}$ of $I$, it seems there is natural way in which the exponents on the variables of these generators change.  Though this seems to be quite difficult to classify.  In every instance, it seems that the row operations will have no affect on whether or not a row in $f_3$ is dependant mod $I$.  There are also no other known examples of ideals where you get non pure power entries from $I$ in $f_3$.  Because of this we offer the following conjecture.

\begin{conjecture}
\label{IRes}
Let $I$ be an $\mathfrak{m}$-primary, Artinian monomial ideal minimally generated by $n$ monomials, then the following hold.
\begin{enumerate}
\item If $f_3$ contains nonzero entries from $I$ these must be pure power generators of $I$.
\item $\mu_R^1=m+n-1-r$ where $r$ is the number of nonzero entries in $f_3$ from $I$.
\end{enumerate}
\end{conjecture}

If Cojecture~\ref{IRes} can be proved then the proof of the following conjecture would be straightforward from the results in this paper.  This would give us a complete classification for the Koszul algebra of Artinian, trivariate monomial ideals.

\begin{conjecture}
\label{CompClass}
Let $I$ be an $\mathfrak{m}$-primary, Artinian monomial ideal minimally generated by $n$ monomials, then we have the following Koszul algebra classifications for $R$:

\begin{enumerate}

\item If $I=\brc{\small x^a,y^b,z^c,x^{a_1}y^{b_1}z^{c_1},\ldots x^{a_{\rho}}y^{b_{\rho}}z^{c_{\rho}}}$ with $a,b,c,a_i,b_i,c_i>0$ for all $1\leq i\leq \rho$ and $\rho>0$, then $R$ is in \textbf{\textup{T}}.
\item If $I=\brc{x^a,y^b,z^c,\{x^{\alpha_{1i}}y^{\alpha_{2i}}z^{\alpha_{3i}}\}_{i=1}^{\rho},\{x^{a_{1i}}y^{a_{2i}}z^{a_{3i}}\}_{i=1}^{\beta}}$ with $\rho\geq2$, satisfying the following,
\begin{enumerate}
\item for exactly one $i\in\{1,2,3\}$ we have that $\alpha_{ij}=\alpha_{ik}>0$ for all $k,j\in\{1,\ldots,\rho\},$ and $a_{il}>\alpha_{ij}$ for all $1\leq l \leq \beta$,
\item for $s,j\in\{1,2,3\}$, $s\not=t$ and $s,t\not=i$ we have that $\alpha_{s1},\alpha_{t\rho} = 0$, $\alpha_{sj} < \alpha_{sj+1}$ and $\alpha_{tk} > \alpha_{tk+1}$ for every $j,k\in\{1,\ldots,\rho-1\}.$
\end{enumerate}
Then $R$ is in \textbf{\textup{B}}.
\item If $I=\brc{x^a,y^b,z^c,x^{\alpha_1}y^{\beta_1},x^{\alpha_2}z^{\gamma_1},y^{\beta_2}z^{\gamma_2},x^{a_1}y^{b_1}z^{c_1},\ldots, x^{a_{\rho}}y^{b_{\rho}}z^{c_{\rho}}}$,\\ $a,b,c,\alpha_i,\beta_i,\gamma_i>0$, then $R$ is in \textbf{\textup{H}}$(0,0)$ and is Golod if and only if there are no $\sigma_{ij}\in S_2$ with only nonzero entries from $I$.
\item Otherwise $R$ is in \textbf{\textup{H}}$(p,q)$ with $p+q\geq 1$.
\end{enumerate}

\end{conjecture}

\end{document}